\documentclass{amsart}

\usepackage{amscd,amssymb}

\def\Z{{\mathbb Z}}
\def\Q{{\mathbb Q}}
\def\R{{\mathbb R}}
\def\C{{\mathbb C}}
\def\P{{\mathbb P}}

\def\A{{\mathcal A}}

\def\J{{\mathcal J}}
\def\L{{\mathcal L}}
\def\M{{\mathcal M}}

\def\O{{\mathcal O}}

\def\cC{{\mathcal C}}

\def\p{{\mathfrak p}}

\def\u{{\mathfrak u}}

\def\G{{\Gamma}}
\def\D{{\Delta}}

\def\d{{\delta}}
\def\w{{\omega}}
\def\bw{{\pmb \w}}

\def\dot{{\bullet}}
\def\bs{\backslash}

\def\bil#1#2{\langle #1,#2 \rangle}

\newcommand\ad{\operatorname{ad}}

\newcommand\Der{\operatorname{Der}}
\newcommand\Jac{\operatorname{Jac}}
\newcommand\Pic{\operatorname{Pic}}
\newcommand\Gr{\operatorname{Gr}}

\newtheorem{theorem}{Theorem}
\newtheorem{lemma}[theorem]{Lemma}
\newtheorem{proposition}[theorem]{Proposition}
\newtheorem{corollary}[theorem]{Corollary}

\theoremstyle{definition}

\theoremstyle{remark}
\newtheorem{remark}[theorem]{Remark}

\begin{document}

\title{Geometric Proofs of Some Results of Morita}

\author{Richard Hain}
\address{Department of Mathematics\\ Duke University\\
Durham, NC 27708-0320}
\email{hain@math.duke.edu}

\author{David Reed}
\address{Department of Mathematics\\ Duke University\\
Durham, NC 27708-0320}
\email{dreed@math.duke.edu}

\date{\today}

\thanks{The first author was supported in part by grants from the National
Science Foundation.}

\subjclass{Primary 14H10; Secondary, 20J05}


\maketitle

\section{Introduction}

In this note we consider three results of Morita: \cite[(1.3)]{morita:jac1},
\cite[(1.7)]{morita:jac2} and \cite[(5.1)]{morita:tani}, which give relations
between certain two dimensional cohomology classes of various moduli spaces of
curves. We have reformulated Morita's results in more geometric language to
facilitate their application in our work \cite{hain-reed} on the Arakelov
geometry of $\M_g$ and to advertise them to algebraic geometers. We give
a new, and hopefully more geometric, proof of each result. In the last
section, we give a geometric interpretation of another of Morita's results
\cite[(5.4)]{morita:jac1}.

Denote the mapping class group of a closed orientable surface of genus
$g$ by $\G_g$ and the mapping class group of a pointed, closed orientable
surface by $\G_g^1$.

Denote the moduli space of smooth projective curves over $\C$
(i.e., compact Riemann surfaces of genus $g$)
by $\M_g$ and the universal curve over it by $\cC_g$. We denote the 
moduli space of principally polarized abelian varieties of dimension $g$
by $\A_g$. As is customary, we view each of these as an orbifold, or,
more accurately, a stack in the sense of \cite{mumford}. The orbifold
fundamental groups of $\M_g$, $\cC_g$ and $\A_g$ are isomorphic to
$$
\G_g,\, \G_g^1 \text{ and }Sp_g(\Z)
$$
respectively, (except when $g=1$ when
$\pi_1(\cC_1,*)\cong SL_2(\Z)\ltimes\Z^2$ and
$\G_1=\G_1^1 \cong SL_2(\Z)$);
the isomorphisms are unique up to inner automorphisms. Each of these
moduli spaces is a $K(\pi,1)$ in the orbifold sense, so that for each
of these moduli spaces $X$
$$
H^\dot(X,\Q) \cong H^\dot(\pi_1(X),\Q)
$$
where $\pi_1(X)$ denotes its orbifold fundamental group.
The period map $\M_g \to \A_g$ takes the moduli point $[C]$ of the
$C$ to that of its jacobian $[\Jac C]$. It corresponds to the canonical
homomorphism $\G_g \to Sp_g(\Z)$.

Suppose that $l$ is a positive integer. The moduli space of smooth
projective curves of genus $g$ with a level $l$ structure will be denoted
by $\M_g[l]$ and the universal curve over it by $\cC_g[l]$. Definitions of
these moduli spaces and the corresponding mapping class groups, $\G_g[l]$
and $\G_g^1[l]$, can be found in \cite[\S3,\S7.4]{hain-looijenga}.
Note that for all $l$,
$$
H^\dot(\M_g,\Q) = H^\dot(\M_g[l],\Q)^{Sp_g(\Z/l\Z)}
\text{ and }
H^\dot(\cC_g,\Q) = H^\dot(\cC_g[l],\Q)^{Sp_g(\Z/l\Z)}.
$$
In particular, the cohomology of moduli space with no level is always a
subspace of the cohomology of the moduli space with a level structure.

The {\it Hodge bundle} is the rank $g$ vector bundle over $\A_g$ whose
fiber over the moduli point $[A]$ of the abelian variety $A$ is
$H^{1,0}(A)$. Its first Chern class is denoted $\lambda$, as is its
pullback to $\M_g$ along the period map. It is known that for each
$g\ge 1$, $\lambda$ generates $H^2(\A_g)$ and $H^2(\M_g)$. (In fact,
in both cases, the determinant $\L$ of the Hodge bundle generates the
Picard group of the corresponding moduli functor, \cite{arb-corn}.) 

Let $\w$ be the first Chern class of the relative {\it cotangent} bundle $\bw$
of the projection $\cC_g \to \M_g$. It is well known that $H^2(\cC_g,\Q)$
has basis $\lambda$ and $\w$. The class $\w$ is often denoted by $\psi$
in the physics literature.\footnote{Recall that $\kappa_1 = 12\lambda$.
Morita denotes $\kappa_1$ by $e_1$ and $-\psi$ by $e$.}

One has the universal abelian variety $\J \to \A_g$. This has (orbifold)
fundamental group isomorphic to $Sp(H_\Z)\ltimes H_\Z$, where $H_\Z$
denotes the first homology group of the reference abelian variety. The
splitting is given by the zero section. The pullback of $\J$ along
the period map is the jacobian of the universal curve.

Since the universal abelian
variety is a flat bundle of real tori, the projection is a foliated
bundle. It follows that there is a closed 2-form on $\J$ whose restriction
to each fiber is the invariant form that corresponds to the polarization
(i.e., the invariant bilinear form on $H_1(A)$), that is parallel with
respect to the flat structure, and whose restriction to the zero section
is trivial (cf.\ \cite[\S1]{morita:jac1}). We shall denote the cohomology
class of this form by
$$
\phi \in H^2(\J,\Q) \cong H^2(Sp_g(\Z)\ltimes H_\Z,\Q).
$$
The class $\phi$ is not integral, but $2\phi$ is, \cite[(1.2)]{morita:jac1}.

Denote the canonical divisor of the Riemann surface $C$ by $K_C$.
Define a map $\kappa : \cC_g  \to \J$  from the universal curve to
the universal abelian variety by taking the point $[C,x]$ of $\cC_g$
to
$$
(2g-2)x - K_C \in \Jac C
$$

\begin{theorem}[Morita {\cite[(1.7)]{morita:jac2}}]
\label{riem_const}
For all $g\ge 1$,
$$
\kappa^\ast \phi = 2g(g-1)\w - 6\lambda \in H^2(\cC_g,\Q).
$$
\end{theorem}

A theta characteristic of a compact Riemann surface $C$ is a square
root of its canonical bundle $K_C$; that is, a divisor class $\alpha$ such
that $2\alpha = K_C$. Each theta characteristic $\alpha$ determines a 
divisor $\Theta_\alpha$ in the jacobian of $C$.
If one passes to the level 2 moduli space $\M_g[2]$,
one can consistently choose a theta characteristic for each curve.%
\footnote{Actually, one has such sections defined over the moduli space
of pairs $(C,\alpha)$ where $C$ is a smooth projective curve of genus $g$
and $\alpha$ is a theta characteristic. This moduli space has two components;
one for the even theta characteristics, the other for the odd ones. Since
theta characteristics correspond the $\Z/2$ quadratic forms on $H_1(C,\Z/2)$
associated to the intersection pairing, all theta characteristics will be
defined over $\M_g[2]$, the moduli space of compact Riemann surfaces
with a full level 2 structure. We prefer to work with this space which
dominates the two moduli spaces of curves with theta characteristic.
It is known by the results of Harer
\cite{harer:spin} and Foisy \cite{foisy} that the second cohomology
with $\Q$ coefficients of $\cC_g[2]$ and the spin moduli spaces of 
pointed curves are isomorphic to $H^2(\cC_g,\Q)$, but we do not 
need such deep results in this paper.} In this case we can define a
map
$$
j_\alpha : \cC_g[2] \to \J[2]
$$
by taking $[C,x]$ to $(g-1)x - \alpha \in \Jac C$. Observe that
$\kappa = 2j_\alpha$.

For each such $\alpha$, there is a theta divisor $\Theta_\alpha$
in $\J[2]$, the universal abelian variety over $\A_g[2]$, the moduli
space of principally polarized abelian varieties with a level 2
structure. Any two such
theta divisors differ by a point of order 2 of $\J[2]$. The rational
homology class of $\Theta_\alpha$ is therefore independent of $\alpha$
and is the pullback of a class in $H^2(\J[2],\Q)$ that we denote by
$\theta$.

The following result is presumably well known; a proof is given in
Section~\ref{proof_theta}.

\begin{proposition}
\label{twist}
For all $g\ge 1$,
$$
\phi = \theta - \lambda/2 \in H^2(\J[2],\Q).
$$
\end{proposition}

Since $\phi$ and $\lambda$ both lie in the subspace $H^2(\J,\Q)$, it follows
that $\theta \in H^2(\J,\Q)$.

We shall prove the following result which is easily seen, using the
proposition, to be equivalent to Theorem~\ref{riem_const}.

\begin{theorem}
\label{theta}
For all $g\ge 1$, 
$$
j_\alpha^\ast \theta = \binom{g}{2}\w - \lambda \in H^2(\cC_g[2],\Q).
$$
\end{theorem}

Denote the square 
$$
\cC_g \times_{\M_g} \cC_g
$$
of the universal curve
by $\cC^2_g$. It is the moduli space of triples $(C,x,y)$ where $x$ and
$y$ are arbitrary points of $C$, a compact Riemann surface. There is a
canonical projection $\cC^2_g \to \M_g$. We then have a commutative
square
$$
\begin{CD}
\cC_g^2 @>\d>> \J \cr
@VVV @VVV \cr
\M_g @>>> \A_g \cr
\end{CD}
$$
where $\d$ is the difference map which is defined by
$$
\d : [C,x,y] \mapsto [x] - [y] \in \Jac C.
$$

The diagonal copy of $\cC_g$ in $\cC_g^2$ is a divisor and thus has
a class $\D$ in $H^2(\cC_g^2,\Q)$. For $j=1,2$, denote the first Chern
class of the relative cotangent bundle of the $j$th projection
$p_j : \cC_g^2 \to \cC_g$ by $\psi_j$. (That is, $\psi_j= p_j^\ast(\w)$.)

\begin{theorem}[Morita {\cite[(1.3)]{morita:jac1}}]
\label{difference}
For all $g\ge 1$,
$$
\delta^*(\phi) = \D + (\psi_1 + \psi_2)/2 \in H^2(\cC_g^2,\Q).
$$
\end{theorem}

Combining this with Proposition~\ref{twist} we obtain:

\begin{corollary}
\label{theta2}
For all $g\ge 1$,
$$
\delta^\ast(\theta) =  \D + (\lambda + \psi_1 + \psi_2)/2
\in H^2(\cC_g^2[2],\Q).
$$
\end{corollary}

In order to state the third theorem, we need a generalization of the
construction of $\J$, the universal abelian variety. Suppose that $W_\Z$
is an $Sp_g(\Z)$ module and that the $Sp_g(\Z)$ action on
$W_\R:=W_\Z\otimes \R$ extends to an action of the Lie group $Sp_g(\R)$.
One can then form the corresponding flat (orbifold) bundle $\J(W)$ of
tori over $\A_g$. There are three ways to view this:
\begin{enumerate}
\item $\J(W)$ is
$$
(Sp_g(\Z) \ltimes W_\Z)\bs (Sp_g(\R)\ltimes W_\R)/U(g)
$$
which maps to $\A_g = Sp_g(\Z)\bs Sp_g(\R)/U(g)$.
\item $\J(W)$ is the quotient of the flat bundle over $\A_g$ associated
to $W_\R$ by the local system over $\A_g$ corresponding to $W_\Z$. 
\item If $W_\Z$ underlies a variation of Hodge structure of odd weight
over $\A_g$, then $\J(W)$ can be identified with the corresponding family
of Griffiths intermediate jacobians over $\A_g$.
\end{enumerate}

There are three cases of interest to us:
\begin{enumerate}
\item $W_\Z = H_\Z$ where $\J(W) = \J$, the universal abelian variety.
\item $W_\Z = \Lambda^3 H_\Z$; this is a variation of Hodge structure and
$\J(\Lambda^3 H)$ is the bundle of intermediate jacobians associated with
the local system over $\A_g$ with fiber $H_3(A,\Z)$ over $[A]$.
\item $W_\Z = \Lambda^3 H_\Z/H_\Z$ where $H_\Z$ is imbedded in $\Lambda^3
H_\Z$ by wedging with the polarization $\zeta \in \Lambda^2 H_\Z$.
This is a variation of Hodge structure and $\J(\Lambda^3 H/H)$ is the
bundle of intermediate jacobians associated with the primitive part of
the local system over $\A_g$ whose fiber over $[A]$ is the `primitive
part' of $H_3(A)$.
\end{enumerate}

Suppose now that $W_\Z$ has an $Sp_g(\Z)$-invariant skew symmetric bilinear
form $q$ (as each of the examples above does.) This corresponds to an
$Sp_g(\Z)$ invariant cohomology class $q \in H^2(W_\R/W_\Z,\Z)$.
The projection $\J(W) \to \A_g$
is a foliated bundle of tori. Consequently, there is a closed 2-form on
$\J(W)$ whose restriction to each fiber is the invariant form that
represents the class of $q$, that is parallel with respect to the flat
structure, and whose restriction to the zero section is trivial (cf.\
\cite[\S1]{morita:jac1}). We shall denote the cohomology class of this
form in $H^2(\J(W),\Q)$ by $\phi(W)$. The class $2\phi(W)$ is always
integral.

There are lifts
$
\nu : \M_g \to \J(\Lambda^3 H/H) \text{ and } \mu : \cC_g \to \J(\Lambda^3 H)
$
of the period maps
$$
\M_g \to \A_g \text{ and } \cC_g \to \A_g
$$
such that the diagram
$$
\begin{CD}
\cC_g @>\mu>> \J(\Lambda^3 H) \cr
@VVV		@VVV \cr
\M_g @>\nu>> \J(\Lambda^3 H/H) \cr
\end{CD}
$$
commutes and the composite of $\mu$ with the map 
$\J(\Lambda^3 H) \to \J(H)$ induced by the contraction $c:\Lambda^3 H \to H$
defined in Section~\ref{lin_alg} is the map $\kappa : \cC_g \to \J(H)$ of
Theorem~\ref{riem_const}.
The map $\mu$ is the normal function of the algebraic cycle in the universal
jacobian $\J \to \M_g$ whose fiber over the point $[C,x]$ of $\cC_g$ is the
algebraic cycle $C_x - C_x^-$ in $\Jac C$. More details can be found in 
\cite[\S 6]{hain:comp}, \cite{hain-looijenga} and \cite{pulte}.
The period maps $\nu$ and $\mu$ induce homomorphisms
$$
\G_g \to Sp_g(\Z)\ltimes \left(\Lambda^3 H_\Z/H_\Z\right)
\text{ and }
\G_g^1 \to Sp_g(\Z) \ltimes \Lambda^3 H_\Z
$$
on fundamental groups. In both cases, the restriction of
this homomorphism to the Torelli group is twice the Johnson homomorphism
(\cite[(6.3)]{hain:comp}, see also \cite{hain:normal}).\footnote{The
homomorphism induced by $\mu$ was first constructed by Morita in
\cite{morita:lift}. The equality of Morita's homomorphism and the one
induced by the normal function follows from results in \cite{hain:normal}.}

Each of $H$, $\Lambda^3 H$ and $\Lambda^3 H/H$ has an $Sp_g(\Z)$ invariant 
skew symmetric bilinear form:
\begin{enumerate}
\item The form on $H$ is the intersection pairing, $(x,y)$.
\item The form on $L:=\Lambda^3 H$ is defined by
$$
\bil{x_1\wedge x_2 \wedge x_3}{ y_1\wedge y_2\wedge y_3} = \det(x_i,y_j).
$$
\item The form on $L/H=\Lambda^3 H/ H$ is constructed in
Section~\ref{lin_alg}.
\end{enumerate}
Each of these forms is primitive in the sense that it cannot be divided by
an integer and still be integer valued.

Set
$$
\phi_H = \kappa^\ast \phi(H),\, \phi_L = \mu^\ast \phi(L)
\text{ and } \phi_{L/H} = \nu^\ast \phi(L/H).
$$
Note that $\phi(H)$ equals the class $\phi$ defined previously and that
$\phi_{L/H} = 0$ when $g=1,2$.

The third result of Morita we wish to prove is:

\begin{theorem}[Morita {\cite[(5.1)]{morita:tani}}]
\footnote{One should note that Morita's class $z_1$ equals $2\phi_H$ and
his class $z_2$ equals $3\phi_L$. We will prove this in
Section~\ref{absolute}.}
\label{chern}
For all $g\ge 1$,
$$
\w = \frac{1}{2g(2g+1)}(2 \phi_H + 3 \phi_L)
\in H^2(\cC_g,\Q).
$$
\end{theorem}

By some easy computations in linear algebra (cf.\ Corollary~\ref{reln})
we have
$$
(g-1)\phi_L = \phi_H + \phi_{L/H}.
$$

Combining this with Theorem~\ref{riem_const}, one arrives at the following
equivalent result which we shall prove in Section~\ref{proof_absolute}:

\begin{theorem}[Morita {\cite[(5.8)]{morita:tani}}]
\label{absolute}
For all $g\ge 1$, 
$$
2\phi_{L/H} = (8g+4) \lambda \in H^2(\M_g,\Z).
$$
\end{theorem}

\begin{remark}
\label{tf}
\begin{enumerate}
\item Since $H_1(\G_g^n,\Z)$ is torsion free when $g\ge 3$, there is no
torsion in $H^2(\G_g^n,\Z)$ when $g\ge 3$. So if a relation between integral
cohomology classes holds in $H^2(\G_g^n,\Q)$, it holds in $H^2(\G_g^n,\Z)$.
\item The $(8g+4)\lambda$ appears frequently in identities involving divisor
classes on $\M_g$, such as in the work of Cornalba and Harris
\cite{corn-harris} and Moriwaki \cite{moriwaki}. These papers are related to
Theorem~\ref{absolute} as will be explained in a forthcoming
paper by the first author.
\item Theorem~\ref{absolute} also gives a clean direct proof of 
the non-triviality of the central extension in \cite{hain:comp}.
\end{enumerate}
\end{remark}

\noindent {\it Acknowledgements:} We would like to thank the referee
for making useful suggestions that lead to the simplification and
clarification of the proof of Theorem~\ref{theta}. In particular, the
statement of Theorem~\ref{square} was suggested by the referee.

\section{Weierstrass Points}

Throughout this section, $C$ denotes a smooth projective curve of
genus $g$, where $g>1$. A point $P$ of $C$ is a {\it Weierstrass point} if
there is a non-constant rational function $f : C \to \P^1$ of degree $g$
such that $f^{-1}(\infty) = P$. Equivalently, $P$ is a Weierstrass point
if and only if $h^0(gP) > 1$. The multiplicity of $P$ as a Weierstrass
point is $h^0(gP) - 1$. The Weierstrass point divisor of $C$ is
$$
W_C = \sum_{P\in C} (h^0(\O(gP) - 1)P.
$$
It has degree $g(g-1)(g+1)$.

We shall denote
the group of divisor classes of degree $d$ on $C$ by $\Pic^d C$ and the
$d$th symmetric power of $C$ by $C^{(d)}$. The divisor class map gives a 
canonical morphism $\mu_d:C^{(d)} \to \Pic^d C$ whose image is the locus of
effective divisor classes of degree $d$, which we shall denote by $W_d$. The
map of $C^{(d)}$ to $W_d$ is well known to be of degree 1 when $1\le d \le g$.
It is classical that $\mu_1$ is an imbedding. We shall identify $W_1$ with $C$
via $\mu_1$ and denote $W_{g-1}$ by $\Theta_C$.

Define
$$
\sigma_C : C \times C \to \Pic^{g-1} C
$$
by $\sigma_C(P,Q) = gP - Q$.

\begin{lemma} \label{pullback}
There is a divisor $E_C$ in $C$, with the same support as
$W_C$, such that
$$
\sigma_C^\ast \Theta_C = g\Delta + p_1^\ast E_C,
$$
where $p_1 : C\times C \to C$ is projection onto the first factor.
\end{lemma}

\begin{proof}
First, suppose that $P$ and $Q$ are two distinct points of $C$. If
$gP - Q \in \Theta_C$, then  there exist $P_1,\dots,P_{g-1} \in C$ such
that
$$
gP \equiv Q + P_1 + \dots + P_{g-1}.
$$
Since $P \neq Q$, this implies that $h^0(\O(gP)) > 1$ and that $P$ is a
Weierstrass point. It follows that
$$
\sigma_C^\ast \Theta_C = n \Delta + p_1^\ast E_C
$$
where $n\in \Z$ and $E_C$ is a divisor in $C$ with the same support at $W_C$.

To determine $n$, pick a point $P$ of $C$ that is not a Weierstrass point.
Then, by the computation above, the image of $C$ in $\Pic^{g-1}C$ under the
mapping
$$
\nu_P : Q \mapsto gP - Q
$$
intersects $\Theta_C$ only at $(g-1)P$. The Poincar\'e dual of $\Theta_C$ is
the polarization
$$
a_1\wedge b_1 + \dots + a_g \wedge b_g \in H^2(\Pic^{g-1} C,\Z)
$$
where $a_1,\dots , b_g$ is a symplectic basis of $H^1(C ,\Z)$. Integrating
this over the image of $C$, we see that the degree of the pullback of
$\Theta_C$ to $C$ along $\nu_P$ is $g$. It follows that
$$
\nu_P^\ast \Theta_C = g P
$$ 
and that $n=g$.
\end{proof}

Define
$$
j_C : C \to \Pic^{g-1} C
$$
by $j_C(P) = (g-1)P$. Note that $j$ is just the restriction of $\sigma_C$ to
the diagonal in $C\times C$. Since the normal bundle of the diagonal 
in $C\times C$ is the tangent bundle of $C$, we have:

\begin{corollary}
There is an isomorphism of line bundles
$$
j_C^\ast\O(\Theta_C) \cong \bw_C^{\otimes(-g)}\otimes \O(E_C). \qed
$$
\end{corollary}

\begin{lemma}\label{degree}
We have $\deg E_C = \deg W_C$.
\end{lemma}

\begin{proof}
The degree of $j_C^\ast\Theta_C$ on $C$ is easily seen to be $g(g-1)^2$.
By Lemma~\ref{degree}, the degree of $E_C$ equals the degree of
$\O(\Theta_C)\otimes \bw^{\otimes g}$, which is
$$
2g(g-1) + g(g-1)^2 = g(g-1)(g+1) = \deg W_C.
$$
\end{proof}

We can now apply this to the universal curve $\cC_g$ over $\M_g$. We have
the mapping
$$
\sigma : \cC_g^2 \to \Pic^{g-1}_{\M_g}\cC_g
$$
whose restriction to the fiber over $[C]\in \M_g$ is $\sigma_C$. Denote
the Weierstrass point divisor in $\cC_g$ by $W$ and the projection onto
the first factor by $p_1 : \cC_g^2 \to \cC_g$. The universal theta divisor
in $\Pic^{g-1}_{\M_g}\cC_g$ will be denoted by $\Theta$. It has fiber
$\Theta_C$ over $[C]\in \M_g$.

\begin{theorem}\label{square}
On $\cC_g^2$, we have $\sigma^\ast \Theta = g\Delta + p_1^\ast W$.
\end{theorem}

\begin{proof}
Set $E = \sigma^\ast \Theta - g\Delta$. By Lemma~\ref{pullback}, the
support of $E$ equals that of $W$. By Lemma~\ref{degree}, $E$ and $W$
have the same degree over $\M_g$. Since the Weierstrass points of the
generic curve are distinct, $E$ and $W$ must be equal.
\end{proof}

Restricting $\sigma$ to the diagonal copy of $\cC_g$ in $\cC_g^2$, we
obtain a morphism
$$
j : \cC_g \to \Pic^{g-1}_{\M_g}\cC_g
$$
whose restriction to the fiber over $[C] \in \M_g$ is $j_C$.

\begin{theorem}
\label{wp_theta}
As line bundles over $\cC_g$,
$$
j^\ast \O(\Theta) \cong \bw^{\otimes(-g)}\otimes \O(W).
$$
\end{theorem}

This follows from the previous result and the following lemma.

\begin{lemma}
\label{normal}
The normal bundle of $\cC_g$ imbedded diagonally in $\cC^2_g$ is the relative
tangent bundle $\bw^{-1}$ of $\cC_g \to \M_g$.
\end{lemma}

\begin{proof}
Let $\pi : \cC_g \to \M_g$ denote the projection. The two projections of
$\cC^2_g$ induce a map $T\cC_g^2 \to T\cC_g \oplus T\cC_g$. We then have
an exact sequence
$$
\begin{CD}
0 @>>> T\cC_g @>\Delta>> (T\cC_g \oplus T\cC_g)|_{\text{diag}}
@>f>> T\cC_g @>>> T\M_g @>>> 0
\end{CD}
$$
where $\Delta$ is the diagonal map and $f(v,w) = \pi_\ast(v) - \pi_\ast(w)$.
It follows that the normal bundle is the kernel of $T\cC_g \to T\M_g$,
which is the relative tangent bundle.
\end{proof}

\section{Proof of Theorem \ref{theta}}
\label{proof_theta}

We will prove Theorem~\ref{theta} by comparing the formula for the class of
Weierstrass point locus $W$ derived in the previous section with a second
expression for the class of $W$ derived using Wronskians. This second
expression is due to Arakelov \cite[Lemma 3.3]{arakelov} as was pointed out
to us by Jean-Benoit Bost.

\begin{proposition}[Arakelov]
Over $\cC_g$ we have
$$
c_1(\O(W)) = \binom{g+1}{2} w - \lambda.
$$
\end{proposition}

\begin{proof}
Suppose that $\zeta_1(t), \dots, \zeta_g(t)$ is a local holomorphic framing
of the Hodge bundle $\pi_\ast \Omega^1_{\cC_g/\M_g}$ over $\M_g$. Then 
$$
W(\zeta_1,\dots,\zeta_g)\otimes (\zeta_1\wedge\dots\wedge\zeta_g)^{-1}
$$
is a section of $\bw^{\otimes\binom{g+1}{2}}\otimes \L^{-1}$ that is independent
of the choice of the framing, and therefore extends to a global section over
$\cC_g$. Here
$$
W(\zeta_1,\dots,\zeta_g) :=
\left| \begin{matrix}
f_1(t,z) & \dots & f_g(t,z) \cr
\vdots & \ddots & \vdots \cr
\frac{\partial^{g-1}}{\partial z^{g-1}}f_1(t,z) & \dots
& \frac{\partial^{g-1}}{\partial z^{g-1}}f_g(t,z) \cr
\end{matrix}\right|
dz^{\binom{g+1}{2}}
$$
is the Wronskian, where  $\zeta_j(t) = f_j(t,z) dz$ locally. It is a 
section of $\bw^{\otimes\binom{g+1}{2}}$. The result follows as $W$ is the
divisor of this section by classical Riemann surface theory.
\end{proof}

\begin{proof}[Proof of Theorem~\ref{theta}]
Choose a theta characteristic $\alpha$ defined over $\cC_g[2]$. This gives
an isomorphism
$$
\begin{CD}
\Pic^{g-1}_{\M_g} \cC_g @>{\alpha}>>  \Jac_{/\M_g} \cC_g
\end{CD}
$$
The composition of $j$ with this mapping is the mapping
$$
j_\alpha : \cC_g \to \J = \Jac_{/\M_g} \cC_g
$$
in the introduction.
Denote the image of $\Theta$ in the universal jacobian by $\Theta_\alpha$.
It follows from Theorem~\ref{wp_theta} that
$$
[W] = j_\alpha^\ast\theta + g w \in H^2(\cC_g[2],\Q).
$$
On the other hand, Arakelov's computation implies that
$$
[W] = \binom{g+1}{2} w - \lambda.
$$
Equating these two expressions, we see that
$$
j_\alpha^\ast\theta = \binom{g}{2} w - \lambda.
$$
\end{proof}

\begin{proof}[Proof of Proposition~\ref{twist}]
First note that
$$
H^\dot(\J[2],\Q) \cong H^\dot(Sp_g(\Z)[2]\ltimes H_\Z,\Q)
$$
By a theorem of Raghunathan \cite{raghunathan}, $H^1(Sp_g[2],H_\Q)$
vanishes and, by a theorem of Borel \cite{borel}, $H^2(Sp_g(\Z)[2],\Q)$
has rank one when $g\ge 2$.
Since $\phi$ is not zero in $H^2(\J[2],\Q)$, it follows, by looking at the
Hochschild-Serre spectral sequence of the group extension
$$
0 \to H_\Z \to Sp_g(\Z)[2]\ltimes H_\Z \to Sp_g(\Z)[2] \to 1,
$$
that $H^2(\J[2],\Q)$ is two dimensional and spanned by $\lambda$ and $\phi$.
By restricting the class $\theta$  to any fiber of $\J[2]$, we see that
$$
\theta = \phi + c\lambda.
$$
To determine the rational number $c$, we restrict $\theta$ to the zero section
of $\J[2]$. If we choose $\alpha$ to be even (i.e., $\vartheta_\alpha$ is an
even theta function), then we see that the restriction of $\theta$ to the
zero section is the divisor corresponding to the theta null corresponding
to $\vartheta_\alpha$. Since the theta nulls are sections of a square root
of the determinant $\L$ of the Hodge bundle, as follows from the
transformation law for theta nulls, we see
that $c= 1/2$. The proof in the case when $g=1$ is simpler and left to the
reader.
\end{proof}

\section{Proof of Theorem~\ref{difference}}
\label{proof_difference}

This result is not nearly as deep as the previous result. It follows from
standard computations when $g\ge 3$ that $H^2(\cC_g,\Q)$ is four dimensional
and has basis $\lambda$, $\psi_1$, $\psi_2$ and $\Delta$, the dual of the
diagonally imbedded copy of $\cC_g$. It follows that 
$$
\delta^\ast(\phi) = a\lambda + b\psi_1 + c\psi_2 + d \Delta
$$
where the coefficients are rational numbers. It follows from
Proposition~\ref{normal} that the restriction of $\Delta$ to the diagonal is
$-\psi$. So the restriction of $\delta^\ast(\phi)$ to the diagonal is
$$
a \lambda + (b+c-d)\psi \in H^2(\cC_g,\Q).
$$
The image of the restriction of $\delta : \cC_g^2 \to \J$ to the diagonal
$\cC_g$ has image contained in the zero section of $\J$. Since the
restriction of $\phi$ to the zero section of $\J$ is trivial, it follows that
$a \lambda + (b+c-d)\psi = 0$ in $H^2(\cC_g,\Q)$. So $a=0$ and $b+c-d=0$.
The constants $b$, $c$ and $d$ can now be determined by restricting to any
fiber $C\times C$. In this case, it is an elementary exercise
in algebraic topology to show that under the difference map, $\phi$ pulls
back to $\Delta + (\psi_1 + \psi_2)/2$. The result follows.

\section{Some Linear Algebra}
\label{lin_alg}

As in the introduction, we denote the fundamental representation of
$Sp_g(\Q)$ by $H$ and the standard symplectic form on $H$ by $q_H$.
Denote its third exterior power $\Lambda^3 H$ by $L$ and the natural
symplectic form on it, which was defined in the introduction, by $q_L$.

If $a_1,\dots,a_g,b_1,\dots,b_g$ is a symplectic basis of $H$, then the 
element $\zeta = \sum_{j=1}^g a_j\wedge b_j$ is an $Sp_g(\Q)$ invariant
element of $\Lambda^2 H$ which is independent of the choice of the basis.
Wedging with $\zeta$ induces an $Sp_g(\Q)$ invariant mapping $H\to L$.
One also has the  contraction map $c: L \to H$ which is  defined by
$$
c : x\wedge y \wedge z \mapsto q_H(x,y)z + q_H(y,z)x + q_H(z,x)y.
$$

\begin{proposition}
When $g\ge 2$, the composite $H \to L \to H $ is $g-1$ times the identity.
\qed
\end{proposition}

\begin{corollary}
When $g=2$, $c: L \to H$ is an isomorphism. Consequently, $L/H$ is zero
when $g=1$ and $2$. \qed
\end{corollary}

Denote the image in $L/H$ of $x\wedge y \wedge z \in L$ by
$\overline{x\wedge y \wedge z}$. The projection $p: L \to L/H$ has a
canonical $Sp_g(\Q)$ equivariant splitting $j$ after tensoring with $\Q$.
It is defined by
$$
j(\overline{x\wedge y \wedge z}) = x\wedge y \wedge z -
\zeta \wedge c(x\wedge y \wedge z)/(g-1).
$$
An $Sp_g(\Q)$ invariant symplectic form on $L/H$ is given by
$$
q_{L/H}(u,v) = (g-1) q_L(j(u),j(v)).
$$
One can easily check that this form is integral and primitive. The following
result follows from a straightforward computation.

\begin{proposition}
For all $g\ge 2$ we have $(g-1) q_L = p^\ast q_{L/H} + c^\ast q_H$. \qed
\end{proposition}

\begin{corollary}
\label{reln}
The classes $\phi_L$, $\phi_{H}$ $\phi_{L/H}$ in $H^2(\cC_g,\Q)$, defined in
the introduction, satisfy
$$
(g-1) \phi_L = \phi_{L/H} + \phi_H \qed
$$
\end{corollary}

\section{Proof of Theorem~\ref{absolute}}
\label{proof_absolute}

We begin by dispensing with the cases $g=1$ and $2$. In these cases
$(8g+4)\lambda = 0$ as $H^2(\G_g,\Z)$ is cyclic of order $12$ when $g=1$,
and $10$ when $g=2$. On the other hand, as observed in Section~\ref{lin_alg},
$L/H$ is trivial in these cases too.  So both $(8g+4)\lambda$ and $\phi_{L/H}$
vanish and the result is trivially true.

Now suppose that $g\ge 3$. We shall denote the mapping class group of a compact
oriented surface $S$ of genus $g$ with $r$ boundary components and $n$ distinct
marked points (not lying on the boundary of $S$) by $\G_{g,r}^n$. As usual, $r$
and $n$ are omitted when they are zero. Here we consider only  $\G_g$, $\G_g^1$
and $\G_{g,1}$. The natural homomorphisms
$$
\G_{g,1} \to \G_g^1 \to \G_g
$$
induce homomorphisms
$$
H^2(\G_g,\Z) \hookrightarrow  H^2(\G_g^1,\Z)
\twoheadrightarrow H^2(\G_{g,1},\Z).
$$
Harer \cite{harer:h3} has proved that the left and right hand groups have
rank one and are generated by $\lambda$ for all $g\ge 3$. The middle group
has rank two and is generated by $\lambda$ and $\psi$ for all $g \ge 3$ as
can be seen by applying the Gysin sequence to the group extension
$$
0 \to \Z \to \G_{g,1} \to \G_g^1 \to 1.
$$
Since these groups are torsion free when $g\ge 3$ (cf.\ Remark~\ref{tf}),
and since $2\phi_{L/H}$ is integral, it suffices to show that
$\phi_{L/H} = (4g+2)\lambda$ in $H^2(\G_g^1,\Q)$.

We know that
$$
\phi_{L/H} = x\lambda \in H^2(\G_g,\Q)
$$
for some $x \in \Q$. To determine $x$, we work in $H^2(\G_g^1,\Q)$.
By Theorem~\ref{riem_const},
$$
\phi_H = 2g(1-g)\psi - 6\lambda \in H^2(\G_g^1,\Q),
$$
so that $6\phi_{L/H} - x \phi_H$ is a multiple of $\psi$.

To compute the constant $x$ we convert the problem into linear algebra.
Denote the Lie algebra of the prounipotent radical of the completion
of $\G_{g,r}^n$ with respect to the standard representation
$\G_{g,r}^n \to Sp_g(\Q)$ by $\u_{g,r}^n$. (See \cite{hain:torelli} for
definitions.) There is a natural homomorphism
$$
H^2(\u_{g,r}^n)^{Sp_g} \to H^2(\G_{g,r}^n,\Q)
$$
which is an isomorphism when $g \ge 3$.\footnote{One can see this by
first reducing to the case where $r=n=0$ using the exact sequence
\cite[(3.6)]{hain:torelli}, the result \cite[(8.2)]{hain:torelli},
and a  spectral sequence argument. The case $r=n=0$ can then be proved
using \cite{harer:h3} for the computation of $H^2(\G_g,\Q)$,
\cite[(7.1)]{hain:torelli} for the injectivity of the homomorphism, and
\cite[(11.1)]{hain:torelli} to see that $H^2(\u_g)^{Sp_g}$ is one
dimensional.} 

Since the kernel of $H^2(\G_g^1,\Q) \to H^2(\G_{g,1},\Q)$ is spanned by
$\psi$, it follows that
$$
\Q\,\psi = \ker\{H^2(\u_g^1)^{Sp_g} \to H^2(\u_{g,1})^{Sp_g} \}.
$$
To determine the combinations of $\psi_H$ and $\psi_{L/H}$ that lie in the
kernel of this map we use the short exact sequence
$$
\begin{CD}
0 @>>> \Gr^W_{-2} H_2(\u_{g,1}) @>{\text{cup}^\ast}>> \Lambda^2 H_1(\u_{g,1})
@>\beta>> \Gr^W_{-2}(\u_{g,1}) @>>> 0.
\end{CD}
$$
The first map is the dual of the cup product and $\beta$ is induced by
the bracket. Taking invariants, we obtain the exact sequence
$$
\begin{CD}
0 \longrightarrow \Gr^W_{-2}H_2(\u_{g,1})^{Sp_g} @>{\text{cup}^\ast}>>
\Lambda^2 H_1(\u_{g,1})^{Sp_g}
@>\beta>> \Gr^W_{-2}(\u_{g,1})^{Sp_g} \longrightarrow 0
\end{CD}
$$
Denote the fundamental representation of $Sp_g(\Q)$ by $H$. It follows from
Johnson's Theorem \cite{johnson:h1} that the Johnson homomorphism induces
an $Sp_g(\Q)$ equivariant isomorphism
$$
H_1(\u_{g,1}) \cong \Lambda^3 H.
$$
Let $a_1,\dots,a_g,b_1,\dots,b_g$ be a symplectic basis of $H$
and $A_1,\dots, A_d,B_1,\dots,B_d$ be a symplectic basis of the kernel
of the contraction $c: \Lambda^3 H \to H$ with respect to the symplectic
form $q_L$ on $\Lambda^3 H$, where
$$
2d = \dim \Lambda^3H/H = \binom{2g}{3} - 2g = \frac{2g(2g+1)(g-2)}{3}.
$$
Note that this is where the $8g+4=4(2g+1)$ in the statement of the theorem
will come from. Set
$$
\zeta = \sum_{j=1}^g a_j \wedge b_j \in \left(\Lambda^2 H\right)^{Sp_g}
$$
and
$$
f_H = \sum_{j=1}^g a_j\wedge \zeta \bigwedge b_j \wedge \zeta \text{ and }
f_{L/H} = \sum_{j=1}^d A_j \bigwedge B_j
$$
which are elements of $\left(\Lambda^2 \Lambda^3 H\right)^{Sp_g}$ and
are viewed as elements of $\Lambda^2 H_1(\u_{g,1})^{Sp_g}$ via the
isomorphism above.

Denote the Malcev Lie algebra of the fundamental group of a surface with one
boundary component by $\p_{g,1}$. (It is isomorphic to the free pronilpotent
Lie algebra generated by $H$.) View $\zeta$ as the element $\sum_j [a_j,b_j]$ of
$\Gr^W_{-2} \p_{g,1}$. The canonical map
$$
\alpha : \Gr^W_{-2} \u_{g,1} \to \Gr^W_{-2} \Der \p_{g,1}
$$
is injective.

It follows from a direct computation using \cite[\S 11]{hain:torelli}
or the fact that the Euler class of a surface of genus $g$ is $2-2g$ that
$$
\alpha \circ \beta(f_H) = (2-2g)\ad(\zeta),
$$
and from \cite[(11.4)]{hain:torelli} that
$$
\alpha \circ \beta(f_{L/H})
= -6\,\sum_{j=1}^d\frac{q_{L/H}(A_j,B_j)}{g(2g+1)}\ad(\zeta)
= \frac{-6d}{g(2g+1)}\ad(\zeta)
= -4(g-2)\ad(\zeta)
$$
Consequently
$$
\alpha\circ \beta\left(2(g-2) f_H - (g-1) f_{L/H} \right) = 0
\text{ in } \Gr^W_{-2}\Der \p_{g,1}
$$
and therefore $2(g-2) f_H - (g-1) f_{L/H}$ spans the one dimensional
vector space $H_2(\u_{g,1})^{Sp_g}$. Consequently
\begin{equation}\label{constraint}
\langle 6\phi_{L/H} + x\phi_H, 2(g-2) f_H - (g-1) f_{L/H} \rangle = 0.
\end{equation}

Now $\phi_L\in H^2(\u_g^1)^{Sp_g}$ takes the value
$$
(c^\ast q_H)(f_H) = (g-1)^2 q_H(\zeta) = g(g-1)^2
$$
on $f_H$, while $\phi_{L/H}\in H^2(\u_g^1)^{Sp_g}$ takes the value 
$$
(p^\ast q_{L/H})(f_{L/H}) = (g-1)q_L(f_{L/H}) = d(g-1)
$$
on $f_{L/H}$.
Since $\phi_H(f_{L/H}) = \phi_{L/H}(f_H) = 0$ (by representation
theory, for example), (\ref{constraint}) becomes
$$
g(g-2)(g-1)^2x - 6d(g-1)^2 = 0
$$
from which it follows that $x=4g+2$.

\begin{remark}
We give a brief explanation of why Morita's classes $z_1$ and $z_2$ defined
in \cite[p.~173]{morita:tani} are related to ours by
$$
z_1 = 2\phi_H \text{ and } z_2 = 3\phi_L.
$$
The first point is that his forms $C_1$ and $C_1$ are related to ours by
$$
C_1 = 4 c^\ast q(H) \text{ and } C_2 = 6 q(L).
$$
This can be seen by direct computation.

The second point is that Morita uses these quadratic forms to construct
extensions of $\frac{1}{2}\Lambda^3 H$ by $\Z$, whereas we are considering
extensions of $\Lambda^3 H$ by $\Z$. The isomorphism
$$
(\text{multiplication by }2) : \frac{1}{2}\Lambda^3 H \to \Lambda^3 H
$$
multiplies quadratic forms by 4. So if we rewrite Morita's proof using the
group
$$
Sp_g(\Z) \ltimes \Lambda^3 H \text{ instead of }
Sp_g(\Z) \ltimes \frac{1}{2} \Lambda^3 H
$$
we need to replace $C_1$ by $C_1/4$ and $C_2$ by $C_2/4$.

The final point is that Morita's construction of an extension of
$\Lambda^3 H$ by $\Z$ from a skew symmetric form $q$ (given at the top
of page~174 of \cite{morita:tani}) yields one whose Chern class
is twice the cohomology class in $H^2(\Lambda^3 H,\Z)$ corresponding to $q$.
This may be seen  by computing commutators in the corresponding Heisenberg
groups.

Putting all this together, we have $z_1 = 2  (4\phi_H/4) = 2 \phi_H$ and
$z_2 = 2 (6 \phi_L/4) = 3 \phi_L$.
\end{remark}

\section{Remarks on \cite[(5.4)]{morita:jac1}}

This result of Morita states that if $g\ge 2$, then
$$
H_1(\G_g,H_\Z) \cong \Z/(2g-2)\Z.
$$
Since $H_\Z$ has an $Sp_g(\Z)$ invariant unimodular form, it is isomorphic
to its dual as a $\G_g$ module. So, by the universal coefficient theorem,
$$
\text{torsion subgroup of }H^2(\G_g,H_\Z) \cong \Z/(2g-2)\Z.
$$

\begin{proposition}
For all $g\ge 1$, $H^2(\G_g,H_\Q)$ vanishes.
\end{proposition}

\begin{proof}
Harer \cite{harer:h3} has shown that the natural homomorphism
$\G_{g,1} \to \G_g$ induces an isomorphism on $H^2$ with $\Q$
coefficients for all $g\ge 2$. One can now use the Gysin sequence
$$
0 \to H^0(\G_{g,1},\Q) \to H^2(\G_g^1,\Q) \to H^2(\G_{g,1},\Q)
\to H^1(\G_g^1,\Q) \dots
$$
and the fact that the second Betti number of $\G_g$ is 0 in genus 2
and 1 when $g\ge 3$ \cite{harer:h3} to show that
$$
\dim H^2(\G_g^1,\Q) = 1 + \dim H^2(\G_g,\Q) = 
\begin{cases}
1 & g=2; \cr
2 & g\ge 3. \cr
\end{cases}
$$
Now consider the Hochschild-Serre spectral sequence of the group
extension 
$$
 1 \to \pi_1(\text{reference surface},*) \to \G_g^1 \to \G_g \to 1.
$$
By passing to subgroups of level $\ge 4$ and applying Deligne's 
degeneration result \cite{deligne}, we see that this spectral 
sequence (for cohomology with $\Q$ coefficients) degenerates at
the $E_2$ term. This implies that
$$
\dim H^2(\G_g^1,\Q) = 1 + \dim H^2(\G_g,\Q) + \dim H^1(\G_g,H_\Q).
$$
It follows that $H^1(\G_g,H_\Q)$ vanishes for all $g\ge 2$. The
result is easily proved when $g=1$ using the ``center kills trick''
--- $-I$,  being central in $\G_1 = SL_2(\Z)$ acts
trivially on $H^\dot(\G_1,H_\Q)$ on the one hand, and as $-I$ on the
other as it acts this way in $H$, which forces $H^\dot(SL_2(\Z),H_\Q)$
to be trivial.
\end{proof}

Combining this with Morita's computation, we obtain:

\begin{corollary}
For all $g \ge 2$, we have $H^2(\G_g,H_\Z) \cong \Z/(2g-2)\Z$. \qed
\end{corollary}

This result has a natural geometric interpretation: $H^2(G,A)$, where $G$
is a group and $A$ a $G$-module, classifies extensions of $G$ by $A$; the
identity being the split extension $G\ltimes A$. So Morita's result says
that when $g\ge 2$, the set of equivalence classes of central extensions
of $\G_g$ by $H_\Z$ is a cyclic group of order $2g-2$. It is easy to
realize these extensions geometrically.

For each $d\in \Z$ one has the degree $d$ part $\Pic_{\M_g}^d \cC_g \to \M_g$
of the  relative Picard group of the universal curve. Since the base and fiber
of this map are Eilenberg-MacLane spaces $K(\pi,1)$ (in the orbifold sense),
the total space is too, and we have a short exact sequence of fundamental
groups
$$
0 \to H_\Z \to \pi_1(\Pic_{\M_g}^d \cC_g,*) \to \G_g \to 1
$$
as, for each curve $C$, there is a canonical isomorphism
$H_1(\Pic_{\M_g}^d \cC_g) \cong H_1(C)$. Taking $d$ to the class of this
extension in $H^2(\G_g,H_\Z)$, one obtains a function
$$
\epsilon : \Z \to H^2(\G_g,H_\Z).
$$
Using the addition map
$$
\Pic_{\M_g}^d \cC_g \times_{\M_g} \Pic_{\M_g}^e \cC_g \to \Pic_{\M_g}^{d+e} \cC_g,
$$
it is easy to show this function is a group homomorphism. The canonical
bundle gives a section of $\Pic_{\M_g}^{2g-2} \cC_g$. Consequently, this
homomorphism factors through $\Z/(2g-2)\Z$.

\begin{proposition}
If $g\ge 3$, $\epsilon(d) = 0$ if and only if $d$ is divisible by $2g-2$.
\end{proposition}

\begin{proof}
We have already proved that if $d$ is divisible by $2g-2$, then
$\epsilon(d)=0$. Suppose that $0 < d \le 2g-2$ generates $\ker \epsilon$.
Then $d |(2g-2)$. Set $e = (2g-2)/d$. Since $\epsilon(d) = 0$, there is a
smooth section $s$ of $\Pic_{\M_g}^d \cC_g \to \M_g$. Then
$e\cdot s - (\text{canonical bundle})$
is a smooth section of the universal jacobian and therefore determines a
cohomology class in $H^1(\G_g,H_\Z)$.  It follows from Johnson's theorems
that this group is  torsion (see, for example, \cite[(5.2)]{hain:normal}.)
By the universal coefficient theorem, $H^1(\G_g,H_\Z) \cong H_0(\G_g,H_\Z)$,
which is trivial. It follows that the class of the section
$e\cdot s - (\text{canonical bundle})$ is trivial, and therefore that
$s$ is homotopic to a section of $\Pic_{\M_g}^{d} \cC_g$ corresponding
to an $e$ th root of the canonical bundle. But by the solution
to the Francetta Conjecture \cite{arb-corn} (see also
\cite[\S 12]{hain:normal}), this implies that $d$ is divisible by $2g-2$,
and therefore that $d=2g-2$.
\end{proof}

Combining this with Morita's Theorem, we obtain:

\begin{theorem}
If $g\ge 3$, then the map $\Z/(2g-2)\Z \to H^2(\G_g,H_\Z)$ that takes $d$
to the class corresponding to $\Pic_{\M_g}^d \cC_g$ is an isomorphism. \qed
\end{theorem}

Similarly, one can use the solution of the Francetta Conjecture for moduli
spaces of  curves with a level \cite[\S 12]{hain:normal} to prove that if
$g\ge 3$ the image of the restriction map
$$
H^2(\G_g,H_\Z) \to H^2(\G_g[l],H_\Z)
$$
is $\Z/(2g-2)\Z$ if $l$ is odd and $\Z/(g-1)\Z$ if $l$ is even and
positive, where $\G_g[l]$ denotes the level $l$ subgroup of $\G_g$.
We do not know how to prove surjectivity when $l>0$.

\end{document}